\documentclass[12pt,oneside,english]{amsart}
\usepackage[latin1]{inputenc}
\usepackage{amssymb}

\makeatletter

\usepackage{babel}

\makeatother
\begin{document}
\baselineskip=17pt

\title[On gold ratio]{Gold ratio and a trigonometric identity}

\author{Vladimir Shevelev}

\address{Departments of Mathematics \\Ben-Gurion University of the
 Negev\\Beer-Sheva 84105, Israel. e-mail:shevelev@bgu.ac.il}


\begin{abstract}
We give two proofs of the identity
$$\sqrt{\frac{\cos\frac{2\pi}{5}} {\cos\frac{\pi}{5}} }+\sqrt{\frac{\cos\frac{\pi}{5}} {\cos\frac{2\pi}{5}} }=\sqrt{5},$$
using and not using the gold ratio.
\end{abstract}

\maketitle
\section{INTRODUCTION }\label{s1}
The \slshape gold ratio \upshape $\varphi=1.61803398...$ is one of the most known and astonishing number in mathematics and culture. Its very interesting story is described in, e.g., [1]. It is defined as the positive root of equation

\begin{equation}\label{1}
x^2-x-1=0,
\end{equation}
such that
\begin{equation}\label{2}
\varphi=\frac {1+\sqrt{5}} {2}.
\end{equation}
The following formulas for the gold ratio are very known:
\begin{equation}\label{3}
 \varphi=2\cos \frac {\pi} {5}=\frac {1} {2} \csc\frac {\pi} {10}=1+2\sin\frac {\pi} {10},...
\end{equation}
In this note we give two proofs of a trigonometric identity
\begin{equation}\label{4}
\sqrt{\frac{\cos\frac{2\pi}{5}} {\cos\frac{\pi}{5}} }+\sqrt{\frac{\cos\frac{\pi}{5}} {\cos\frac{2\pi}{5}} }=\sqrt{5}
 \end{equation}
 using and not using the gold ratio. Note that identity (4) has a Ramanujan type structure (cf. [2,p. 326]).
 \section{FIRST PROOF OF IDENTITY 4 }\label{s2}
 Multiplying (1) on $x^{n-1}, \enskip n\geq0,$ we have
 \begin{equation}\label{5}
 x^{n+1}=x^n+x^{n-1}.
 \end{equation}
 Using consequently (5) for $n=2 $ and  $n=1, $  we find
 $$x^3=x^2+x=2x+1.$$
 Thus\newpage
 $$x(x^2-2)=1,\enskip x=\frac {1} {x^2-1}, \enskip x^2=\frac {x} {x^2-1},\enskip x=\sqrt{\frac {x} {x^2-2}}$$
 and
 \begin{equation}\label{6}
  x=\sqrt{\frac {\frac {x}{2} } {\frac {x^2} {2}-1}}.
  \end{equation}
  According to first formula (3), $\varphi=x_1=2\cos \frac {\pi} {5}$ satisfies (6), and we find
  \begin{equation}\label{7}
 \varphi=\sqrt{\frac{\cos\frac{\pi}{5}} {2\cos^2\frac{\pi}{5}-1} }=\sqrt{\frac{\cos\frac{\pi}{5}} {\cos\frac{2\pi}{5}} }.
  \end{equation}
  Now in order to prove (4) it is sufficient to verify the equality
  \begin{equation}\label{8}
  \varphi+\frac{1} {\varphi}=\sqrt{5}.
  \end{equation}
  Indeed, the equation
  \begin{equation}\label{9}
 y^2-\sqrt{5}y+1=0
 \end{equation}
 has the roots $y_1=\frac{\sqrt{5}+1} {2}=\varphi$ and $y_2=\frac{\sqrt{5}-1} {2}=\frac {2} {\sqrt{5}+1} =\frac {1} {\varphi},$ and, by the Vieta theorem, we obtain (8).$\blacksquare$
  \section{SECOND PROOF OF IDENTITY 4 }\label{s3}
  In the second proof we follow to the scheme of our work [3].
 Note that
\begin{equation}\label{10}
\cos\alpha\cos2\alpha=\frac {\sin2\alpha cos2\alpha} {2\sin\alpha }=\frac  {\sin4\alpha} {4\sin\alpha}.
\end{equation}
Substituting in (10) $\alpha=\pi/5,$ we find
\begin{equation}\label{11}
\cos \frac {\pi}{5}\cos \frac {2\pi} {5}=\frac {\sin\frac{4\pi}{5} }{4\sin\frac{\pi} {5}}=1/4.
\end{equation}
On the other hand, we have
$$\cos \frac {2\pi}{5}-\cos \frac {\pi}{5}=\cos \frac {2\pi}{5}+\cos \frac {4\pi}{5}=\frac {2\sin \frac {\pi}{5}}{2\sin \frac {\pi}{5}}(\cos \frac {2\pi}{5}+\cos \frac {4\pi}{5})=$$
\begin{equation}\label{12}
\frac {1}{2\sin \frac {\pi}{5}}(\sin \frac {3\pi}{5}-\sin \frac {\pi}{5}+\sin\pi-\sin \frac {3\pi}{5})=-\frac {1} {2}.
\end{equation}
Putting
\begin{equation}\label{13}
x=-\cos \frac {\pi}{5},\enskip y=\cos \frac {2\pi}{5},
\end{equation}
we have
\begin{equation}\label{14}
xy=-\frac {1} {4},\enskip x+y=-\frac{1} {2}.
\end{equation}
Denoting
\begin{equation}\label{15}
A=\sqrt{-\frac{x} {y}}+\sqrt{-\frac{y} {x}},
\end{equation}\newpage
we have

$$-\frac{x} {y}-\frac{y} {x}=A^2,$$
or
$$2-A^2=\frac {(x+y)^2-2xy} {xy}=-4(\frac{1} {4}+\frac{1} {2})=-3.$$
Thus, since $A>0,$
$$A=\sqrt{5}.$$
Therefore, from (13) and (15) we find (4). $\blacksquare$

\end{document}